\documentclass[12pt,a4paper]{article}
\usepackage[T1]{fontenc}
\usepackage[english]{babel}
\usepackage[utf8]{inputenc}

\usepackage{amsmath}
\usepackage{amsthm}
\usepackage{amsfonts}
\usepackage{amssymb}
\usepackage{yhmath}
\usepackage{MnSymbol}

\usepackage[lmargin=3cm, rmargin=2.5cm,tmargin=2.20cm,bmargin=2.90cm,foot=1.3cm]{geometry}
\usepackage{indentfirst}
\usepackage{orcidlink}
\usepackage{csquotes}
\usepackage{hyperref}

\usepackage[backend=biber,bibencoding=utf8, style=ieee]{biblatex}
\theoremstyle{plain}
\newtheorem{theorem}{Theorem}[section]
\newtheorem{lemma}[theorem]{Lemma}
\newtheorem{corollary}[theorem]{Corollary}
\theoremstyle{definition}
\newtheorem{remark}[theorem]{Remark}
\newtheorem{example}[theorem]{Example}
\newtheorem{question}[theorem]{Question}

\renewcommand{\geq}{\geqslant}
\renewcommand{\leq}{\leqslant}

\newcommand{\rr}{{\mathbb{R}}}
\renewcommand{\ss}{{\mathbb{S}}}
\newcommand{\pol}{\mathbf{p}}
\newcommand{\ran}{\mathsf{ran}}
\begin{document}
\vspace{-40pt}
\title{\textbf{Simple closed curves contained in~$\varepsilon$-boundaries of planar sets}\footnote{The work was performed as part of research conducted in the Ural Mathematical Center with the financial support of the Ministry of Science and Higher Education of the Russian Federation (Agreement number 075-02-2026-737).}
}

\author{Mikhail Patrakeev\footnote{Krasovskii Institute of Mathematics and Mechanics of UB RAS, 620108, Yekaterinburg, Russia; \textit{e-mail}\textup{:} p17533@gmail.com; \href{https://orcid.org/0000-0001-7654-5208}{orcid.org/0000-0001-7654-5208}}\hspace{2mm}\orcidlink{0000-0001-7654-5208}, Aleksei Volkov\footnote{Krasovskii Institute of Mathematics and Mechanics of UB RAS, 620108, Yekaterinburg, Russia; \textit{e-mail}\textup{:} volkov@imm.uran.ru; \href{https://orcid.org/0009-0001-7337-2695}{orcid.org/0009-0001-7337-2695}}\hspace{2mm}\orcidlink{0009-0001-7337-2695}
}
\date{}

\providecommand{\keywords}[1]
{
  \small	
  \textbf{Keywords:} #1
}
\providecommand{\MSC}[1]
{
  \small	
  \textbf{MSC:} #1
}
\vspace{-30pt}
\maketitle
\vspace{-20pt}
\begin{abstract}\rm
\medskip
The $\varepsilon$-boundary of a set ${A}\subseteq\rr^2$ is the set $\{{p}\in\rr^2:\rho({p},{A})=\varepsilon\}$, where $\rho$ is the Euclidean distance.
We prove that if ${A},{B}\subseteq\rr^2$ are nonempty, connected sets, ${A}$ is bounded, and $0<\varepsilon<\rho({A},{B})$, then the $\varepsilon$-boundary of ${A}$ contains a simple closed curve (aka a Jordan curve) that separates ${A}$ and ${B}$.
This statement follows from the theorem which says that if $\varepsilon>0$ and ${A}\subseteq\rr^2$ is a nonempty, bounded, connected set, then the boundary of each component of $\{{p}\in\rr^2: \rho({p},{A})>\varepsilon\}$ is a simple closed curve.
Another corollary of this theorem is that the $\varepsilon$-boundary of a nonempty, bounded, connected set ${A}\subseteq\rr^2$ contains a simple closed curve bounding the domain that contains the open $\varepsilon$-neighbourhood of ${A}$.
In all these statements the connectivity condition can be significantly weakened.
We also show that, for all $\varepsilon>0$, the $\varepsilon$-boundary of a nonempty, bounded set ${A}\subseteq\rr^2$ contains a simple closed curve.
\end{abstract}

\MSC{54E35; 57K20.}


\section{Introduction}
The $\varepsilon$-boundary of a set ${A}\subseteq\rr^2$ is the set $\{{p}\in\rr^2:\rho({p},{A})=\varepsilon\}$, where $\rho$ is the Euclidean distance.
The separation of planar sets by a simple closed curve lying in the $\varepsilon$-boundary of one of these sets is used in engineering.
But the algorithms for constructing such a curve are either based on heuristics or assume some smoothness of the boundary of the sets~\cite{cad,mathpages}.
Similarly, the line of the boundary of territorial waters is drawn heuristically~\cite{maritime}.
No proof of the existence of a simple closed curve that separates arbitrary continua ${A}$ and ${B}$ and lies in the $\varepsilon$-boundary of ${A}$ appears to be known.

\smallskip
It is known that two disjoint continua in the plane can be separated by some simple closed curve~\cite[Chap.\,10, \textsection\,61.II, Theorem\,5$'$]{kuratowski}.
It is also apparently known that if ${A},{B}\subseteq\rr^2$ are nonempty, connected sets, ${A}$ is bounded, the closure of ${A}$ does not separate the plane, and $0<\varepsilon<\rho({A},{B})$, then the open $\varepsilon$-neighbourhood of ${A}$, $\{{p}\in\rr^2:\rho({p},{A})<\varepsilon\}$, contains a simple closed curve that separates ${A}$ and ${B}$.
At least, if we understand R.\,L. Moore's terminology correctly, this statement is a direct consequence of his theorem~\cite[Theorem~1]{moore}, which states that \textit{if, in a plane $S$, $M$ is a closed point set and $K$ is a bounded maximal connected subset of $M$ which does not separate $S$, then, for every positive number $e$, there exists a simple closed curve which encloses $K$ and contains no point of $M$ and which is such that every point within it is at a distance less than $e$ from some point of $K$.}

We prove that if ${A},{B}\subseteq\rr^2$ are nonempty, connected sets, ${A}$ is bounded, and $0<\varepsilon<\rho({A},{B})$, then the $\varepsilon$-boundary of ${A}$ contains a simple closed curve that separates ${A}$ and ${B}$; see Corollary~\ref{cor:boundary.sep} and Remark~\ref{rmk:chain}(1).
In the proof, we specify an example of such a curve explicitly: the boundary of that component of the set $\{{p}\in\rr^2: \rho({p},{A})>\varepsilon\}$ which contains ${B}$.

\smallskip
It is known that, for all $\varepsilon>0$, the open $\varepsilon$-neighbourhood of a nonempty, bounded, connected set ${A}\subseteq\rr^2$ contains a simple closed curve bounding the domain that contains ${A}$.
This result directly follows from a theorem of L. Zoretti \cite{zoretti} which, according to~\cite[Corollary VI.3.11]{whyburn} (see also~\cite{moore} and~\cite{lubben}), states that \textit{if ${K}$ is a component of a compact set $M$ and $\varepsilon$ is any positive number, then there exists a simple closed curve $J$ which encloses $K$ and is such that $J\cap{M} =\varnothing$, and every point of $J$ is at a distance less than $\varepsilon$ from some point of {K}.}

We prove that, for all $\varepsilon>0$, the $\varepsilon$-boundary of a nonempty, bounded, connected set ${A}\subseteq\rr^2$ contains a simple closed curve bounding the region that contains the open $\varepsilon$-neighbourhood of ${A}$; see Corollary~\ref{cor:pzk.in.ball} and Remark~\ref{rmk:chain}(1).
This result admits the following mechanical interpretation: on such a set ${A}$ one can put a wheel of radius $\varepsilon$ lying in the same plane and ``roll it along (stretches of) the boundary of ${A}$'' in such a way that at every moment the wheel touches ${A}$, its interior does not intersect ${A}$, and the center of the wheel eventually describes a simple closed curve bounding the domain that contains the open $\varepsilon$-neighbourhood of ${A}$.

\smallskip
Sometimes there is a need to maximise the Euclidean distance from a simple closed curve to the sets it separates.
We prove that if ${A},{B}\subseteq\rr^2$ are nonempty, connected sets, ${A}$ is bounded, and $\rho({A},{B})>0$, then the $\rho({A},{B})/2$-boundary of ${A}$ contains a simple closed curve ${C}$ that separates ${A}$ and ${B}$ and such that $\rho({C},{A}\cup{B})=\rho({A},{B})/2$; see Corollary~\ref{cor:pzk.sep.rho/2} and Remark~\ref{rmk:chain}(1).
Thus, among the simple closed curves separating ${A}$ and ${B}$, the curve ${C}$ is a curve maximally distant from ${A}\cup{B}$.

We also show that if among the simple closed curves separating nonempty sets ${A}$ and ${B}$ in $\rr^2$, the curve ${C}$ is a curve maximally distant from ${A}\cup{B}$, then either $\rho({C},{A}\cup{B})=0$ or $\rho({C},{A}\cup{B})=\rho({A},{B})/2$, see Remark~\ref{rmk:pzk.sep.max.dist}.
It is curious that there exist two sets ${A}$ and ${B}$ that can be separated by some simple closed curve, but the supremum of the distances from such curves to ${A}\cup{B}$ is not attained; see Subsection~\ref{sec:ex.rho.ab.chain} for an example.
However, we do not know whether there exist nonempty sets ${A},{B}\subseteq\rr^2$ that can be separated by a simple closed curve, with $\rho({A},{B})=1$, and such that, for every simple closed curve ${C}$ that separates ${A}$ and ${B}$, $\rho({C},{A}\cup{B})=0$ holds; see Question~\ref{qst:max.dist}.

\smallskip
In all the above statements, the connectivity condition can be relaxed to the $\delta$-chainedness condition for some $\delta>0$ (a set ${A}$ is \textit{$\delta$-chained} iff any two points of ${A}$ can be connected by a polygonal chain whose vertices belong to ${A}$ and whose segment lengths are all less than $\delta$).
Specific values of the parameter $\delta$ are given in Corollaries~\ref{cor:pzk.in.ball}, \ref{cor:boundary.sep}, and \ref{cor:pzk.sep.rho/2}.

Most of our results follow from Theorem~\ref{thm}, which says that if $\varepsilon>0$ and ${A}\subseteq\rr^2$ is a nonempty, bounded, $2\varepsilon$-chained set, then the boundary of each component of the set $\{{p}\in\rr^2: \rho({p},{A})>\varepsilon\}$ is a simple closed curve.
Almost none of the conditions in Theorem~\ref{thm} and Corollaries~\ref{cor:pzk.in.ball}--\ref{cor:pzk.sep.rho/2} can be relaxed.
Furthermore, most of these statements cannot be transferred to the 3-dimensional case.
We discuss these issues in Sections~\ref{sec:exam} and~\ref{sec:r3}, respectively.
But we do not know whether it is always possible to separate two disjoint continua ${A},{B}\subseteq\rr^3$ by a two-dimensional manifold lying at distance $\rho({A},{B})/2$ from ${A}\cup{B}$; see Question~\ref{qst:r3}.

\smallskip
Finally, we prove that, for all $\varepsilon>0$, the $\varepsilon$-boundary of a nonempty, bounded set ${A}\subseteq\rr^2$ contains a simple closed curve, see Corollary~\ref{cor:pzk.in.bound}.
This statement complements previously known results about $\varepsilon$-boundaries in the Euclidean spaces.
In particular, M. Brown proved~\cite{brown} that, for all $\varepsilon>0$, the $\varepsilon$-boundary of a compact subset of the plane is contained in the union of a finite number of simple closed curves.
In the same paper he showed that if a nonempty, compact set ${A}\subseteq\rr^{n}$ has a diameter smaller than $\varepsilon$ and contains the origin ${0}$, then its $\varepsilon$-boundary is the $({n}\,{-}\,1)$-sphere.
He also established that, for all but countable number of $\varepsilon$, each component of the $\varepsilon$-boundary of a compact subset of the plane is a point, a simple arc, or a simple closed curve.
R. Gariepy and W.\,D. Pepe, answering M. Brown's question, showed~\cite{GP} that the $\varepsilon$-boundary of a closed subset of the plane is a 1-manifold for almost every $\varepsilon$.
Four years later, S. Ferry proved~\cite{ferry} that, for ${n}$ equal to 2 or 3, the $\varepsilon$-boundary of a set ${A}\subseteq\rr^n$ is an $({n}\,{-}\,1)$-manifold for almost all $\varepsilon$.
He also proved that if ${A}$ is a finite polyhedron in $\rr^{n}$, then its $\varepsilon$-boundary is an $({n}\,{-}\,1)$-manifold for all sufficiently small values of $\varepsilon$.
In the same paper he constructed a set ${B}\subseteq\rr^3$ such that the $\varepsilon$-boundary of ${B}$ has components which are not 2-manifolds for uncountably many $\varepsilon$.
And he constructed a Cantor set in $\rr^4$ whose $\varepsilon$-boundary is not a 3-manifold for any $\varepsilon$ between 0 and 1.
P. Pikuta proved~\cite{pikuta} that the $\varepsilon$-boundary of a compact subset of the plane is a closed absolutely continuous curve for all sufficiently large values of $\varepsilon$.
Recently, J. Rataj and L. Zaj\'i\v{c}ek extended~\cite{RZ} the results from~\cite{GP} and~\cite{ferry} to sufficiently smooth normed linear spaces ${X}$ with $\mathsf{dim} X \in\{2,\,3\}$.
\section{Terminology and notation}
We use terminology from the book~\cite{engelking}.
A \textit{simple closed curve} (also called a \textit{Jordan curve}) is a set homeomorphic to a circle (i.e., a one-dimensional sphere).
In metric spaces, simple closed curves are precisely the images of a circle under continuous injective mappings.
According to the Jordan curve theorem, if ${C}$ is a simple closed curve in the plane $\rr^2$, then its complement $\rr^2\setminus{C}$ has exactly two components, the bounded and the unbounded, which we denote by ${C}^{-}$ and ${C}^{+}$, respectively.
We say that a simple closed curve ${C}\subseteq\rr^2$ \textit{separates} sets ${A}$ and ${B}$ iff the sets ${A}$ and ${B}$ are contained in different components of the subspace $\rr^2\setminus{C}$.
We denote the range of a mapping ${f}$ by $\ran({f})$.
A \textit{path} is a continuous mapping whose domain equals the closed segment $[0,1]$.
A path ${f}$ \textit{connects} points ${u}$ and ${v}$ in a space ${X}$ iff ${f}({0})={u}$, ${f}({1})={v}$, and $\ran({f})\subseteq{X}$.
A space ${X}$ is \textit{pathwise connected} if any two of its points are connected by a path in ${X}$.

We define the distance between two sets in $\rr^{n}$ as the infimum of pairwise Euclidean distances between points of these sets.
The distance between a point ${p}$ and a set ${A}$ is the distance between sets $\{{p}\}$ and ${A}$.
We denote all three distances by the symbol $\rho$.
For $\varepsilon>0$ and ${A}\subseteq\rr^2$, the \textit{open $\varepsilon$-neighbourhood} of the set ${A}$, $\{{p}\in\rr^2:\rho({p},{A})<\varepsilon\}$, is denoted by $\mathsf{O}_\varepsilon({A})$;
similarly the \textit{closed $\varepsilon$-neighbourhood}
$\{{p}\in\rr^2:\rho({p},{A})\leqslant\varepsilon\}$ 
and the \textit{$\varepsilon$-boundary} 
$\{{p}\in\rr^2: \rho({p},{A})=\varepsilon\}$ 
of the set ${A}$ are denoted by $\mathsf{B}_\varepsilon({A})$ and $\mathsf{S}_\varepsilon({A})$, respectively.
If ${p}$ is a point in $\rr^2$, then
$\mathsf{O}_\varepsilon({p}) \coloneq \mathsf{O}_\varepsilon(\{{p}\})$,
$\mathsf{B}_\varepsilon({p}) \coloneq \mathsf{B}_\varepsilon(\{{p}\})$, and
$\mathsf{S}_\varepsilon({p}) \coloneq \mathsf{S}_\varepsilon(\{{p}\})$ are the open and the closed disks and the circle of center ${p}$ and radius $\varepsilon$, respectively.
We denote the closure and boundary of a set ${A}$ in $\rr^2$ by $\bar{A}$ and $\partial{A}$, respectively; we denote the boundary of a set ${B}$ in a space ${X}$ by $\partial_{X}{B}$.

\begin{remark}\label{rmk:balls}
Suppose that $\varepsilon>0$, $\varnothing \neq {A}\subseteq\rr^2$, $\varnothing \neq {B}\subseteq\rr^2$, and ${p}\in\rr^2$.
Then:
\begin{enumerate}
\item $\partial\mathsf{B}_\varepsilon({A}) \subseteq \partial\mathsf{O}_\varepsilon({A}) = \mathsf{S}_\varepsilon({A})$.
\item $\rho({p},{A})=\rho({p},\bar{A}) \quad$ and $\quad \rho({A},{B})=\rho(\bar{A},\bar{B})$.
\item $\mathsf{O}_\varepsilon({A}) = \mathsf{O}_\varepsilon(\bar{A}), \quad \mathsf{B}_\varepsilon({A}) = \mathsf{B}_\varepsilon(\bar{A})$, \quad and \quad $\mathsf{S}_\varepsilon({A}) = \mathsf{S}_\varepsilon(\bar{A})$.\hfill$\qed$%
\end{enumerate}
\end{remark}

For $\varepsilon>0$, we say that two points ${p},{q}\in{A}$ are \textit{$\varepsilon$-chained in ${A}$} iff there exists a finite sequence of points ${r}_0, {r}_1,\ldots,{r}_{n}$ in ${A}$ such that ${r}_0={p}$, ${r}_{n}={q}$, and $\rho({r}_{i},{r}_{{i}+1})<\varepsilon$ for all ${i}<{n}$.
A set ${A}$ is called \textit{$\varepsilon$-chained} iff any two of its points are $\varepsilon$-chained in it~\cite[page 60, Definition 4.15]{nadler}.
For ${p}\in{A}$, the \textit{$\varepsilon$-chained component of a point ${p}$ in a set} ${A}$ is the set
\[
\{ {q}\in{A}:\ {p} \text{ and } {q} \text{ are }\, \varepsilon\text{-chained in }{A}\}.
\]
We say that a set ${B}$ is an $\varepsilon$\textit{-chained component of a set} ${A}$ iff ${B}$ equals the $\varepsilon$-chained component of a point ${p}$ in ${A}$ for some ${p}\in{A}$.

\begin{remark}\label{rmk:chain}
Suppose that $\varnothing\neq{A}\subseteq\rr^2$ and $\varepsilon>0$.
Then:
\begin{enumerate}
\item \cite[Exercise 4.23(a)]{nadler} If ${A}$ is connected, then it is $\varepsilon$-chained.
\item If ${A}$ is $2\varepsilon$-chained, then its open $\varepsilon$-neighbourhood $\mathsf{O}_\varepsilon({A})$ is pathwise connected.
\item If ${A}$ is $\varepsilon$-chained and $\delta>\varepsilon$, then ${A}$ is $\delta$-chained.\hfill$\qed$
\end{enumerate}
\end{remark}

Note that in the second clause of Remark~\ref{rmk:chain} the reverse implication is also true.
Thus, a nonempty subset of the plane is $2\varepsilon$-chained if and only if its open $\varepsilon$-neighbourhood is pathwise connected.
\section{Auxiliary lemmas}
To prove the main theorem, we need the following auxiliary statements.

\begin{lemma} \label{lem:connected}
Suppose that ${E}$ is a connected metric space, ${U}\subseteq{E}$ is open, ${p}\in{U}$, and ${U}\setminus\{{p}\}$ is connected.
Then ${E}\setminus\{{p}\}$ is also connected.
\end{lemma}

\begin{proof}
Suppose on the contrary that $E\setminus\{{p}\}$ equals the union of two disjoint nonempty sets $A$ and $B$ closed in $E\setminus\{{p}\}$.
Then there exist two sets ${C}$ and ${D}$ closed in ${E}$ such that 
\[
{A} = {C} \cap ({E}\setminus\{{p}\}) = {C}\setminus\{{p}\}
\qquad \text{and} \qquad
{B} = {D} \cap ({E}\setminus\{{p}\}) = {D}\setminus\{{p}\}.
\]
Clearly, ${C}={A}$ or ${C}={A}\cup\{{p}\}$, and, similarly, ${D}={B}$ or ${D}={B}\cup\{{p}\}$.

\textit{Case 1.}
${C}={A}$ or ${D}={B}$.
Let, without loss of generality, ${C}={A}$.
Note that ${D}\cup\{{p}\}$ is closed in ${E}$, so ${E}$ equals the union of two disjoint nonempty, closed sets ${C}$ and ${D}\cup\{{p}\}$.
This contradicts the connectedness of ${E}$.

\textit{Case 2.}
${C}={A}\cup\{{p}\}$ and ${D}={B}\cup\{{p}\}$.
Note that in this case the sets ${A}$ and ${B}$ are open in ${E}$ as the complements of closed sets ${D}$ and ${C}$, respectively.

Consider the sets
\[
{A}'\coloneq{A}\cap({U}\setminus\{{p}\})
\qquad\text{and}\qquad
{B}'\coloneq{B}\cap({U}\setminus\{{p}\}).
\]

These sets are disjoint and open in ${U}\setminus\{{p}\}$ (because the sets ${A}$ and ${B}$ are disjoint and open in ${E}$) and ${U}\setminus\{{p}\}$ equals their union.
And since, by assumption, ${U}\setminus\{{p}\}$ is connected, one of these sets is empty.
Let, without loss of generality, ${A}'$ be empty.
Then ${U} \subseteq {B}\cup\{{p}\}$, and therefore ${B}\cup\{{p}\} = {B}\cup{U}$, since ${p}\in{U}$.
Thus, the set ${B}\cup\{{p}\}$ is open in ${E}$ as the union of open sets.
Hence, the space ${E}$ equals the union of two disjoint nonempty, open sets ${A}$ and ${B}\cup\{{p}\}$, which contradicts its connectedness.
\end{proof}

\begin{lemma} \label{lem:complement}
Every open, connected set is a component of the complement of its boundary.
\end{lemma}
\begin{proof}
Let ${U}$ be an open, connected set in a topological space ${X}$ and ${C}$ be the boundary of ${U}$.
We need to show that $U$ is a $\subseteq$-maximal connected set in the subspace ${X}\setminus{C}$.
Consider a nonempty set ${V} \subseteq {X}\setminus({U}\cup{C})$.
It suffices to show that ${U}\cup{V}$ is not connected.
Since $U$ is open in ${X}$, then ${U}$ is also open in ${U}\cup{V}$.
The set ${V}$ is also open in ${U}\cup{V}$ because it equals the trace on ${U}\cup{V}$ of the open set ${X}\setminus({U}\cup{C}) = {X}\setminus{\bar{U}}$.
Thus, ${U}\cup{V}$ equals the union of two disjoint nonempty, open sets.
\end{proof}
\section{The main result}
Recall that a \textit{neighbourhood of a point} is a set whose interior contains the given point.
A space is \textit{locally connected} iff every neighbourhood of every point contains a connected neighbourhood of the same point.
A \textit{continuum} is a connected compact set, and a \textit{semi-continuum} is a space such that any pair of points is contained in some continuum.
A point is a \textit {cut point} of a space iff the complement of this point is not a semi-continuum.

\begin{theorem}\label{thm}
Suppose that $\varepsilon>0$ and ${A}\subseteq\rr^2$ is a nonempty, bounded, $2\varepsilon$-chained set.
Let ${D}$ be a component of $\rr^2\setminus\mathsf{B}_\varepsilon({A})$.
Then the boundary of ${D}$ is a simple closed curve contained in the $\varepsilon$-boundary of ${A}$.
Moreover, we have:
\begin{itemize}
\item if ${D}$ is bounded, then ${D}=(\partial{D})^{-}$\textup{;}
\item if ${D}$ is unbounded, then ${D}=(\partial{D})^{+}$.
\end{itemize}
\end{theorem}

\begin{proof}
Let us show that $\partial{D}\subseteq\mathsf{S}_\varepsilon({A})$.
Let ${q}\in\partial{D}$.
If $\rho({q},{A}) < \varepsilon$, then $\mathsf{O}_\delta({q}) \subseteq \mathsf{B}_\varepsilon({A})$ for any $\delta\in\big(0,\varepsilon-\rho({q},{A})\big)$, which contradicts the fact that ${q}\in\partial{D}$.
If $\rho({q},{A}) > \varepsilon$, then $\mathsf{O}_\delta({q}) \subseteq \rr^2\setminus\mathsf{B}_\varepsilon({A})$ for any $\delta\in\big(0,\rho({q},{A})-\varepsilon\big)$.
Thus, the open neighbourhood $\mathsf{O}_\delta({q})$ is a connected subset of the subspace $\rr^2\setminus\mathsf{B}_\varepsilon({A})$ and intersects the component ${D}$ of this subspace, so it is contained in ${D}$ (a connected set is either disjoint with or contained in a component).
Again we get a contradiction with the fact that ${q}\in\partial{D}$.

Let us add a new point $\pol$ (the pole) to the plane $\rr^2$ so that the new space $\ss \coloneq \rr^2\cup\{\pol\}$ is homeomorphic to the two-dimensional sphere.
Being connected, the set ${D}$ is contained in some component ${E}$ of the subspace $\ss\setminus\mathsf{B}_\varepsilon({A})$.

Let us show that
\[
\pol
\notin
\partial_\ss{E}\ \text{ and } \ {E}\setminus\{\pol\} \text{ is connected}.
\]
Since $\pol\notin\mathsf{B}_\varepsilon({A})$ and $\mathsf{B}_\varepsilon({A})$ is closed in $\ss$, there exists a connected open neighbourhood ${U}$ of the point $\pol$ in $\ss$ such that ${U} \cap \mathsf{B}_\varepsilon({A}) = \varnothing$ and ${U}\setminus\{\pol\}$ is connected.
If the connected subset ${U}$ of the subspace $\ss\setminus\mathsf{B}_\varepsilon({A})$ intersects the component ${E}$ of that subspace, then ${U}\subseteq{E}$.
In this case $\pol\notin\partial_\ss{E}$ and, by virtue of Lemma~\ref{lem:connected}, ${E}\setminus\{\pol\}$ is connected.
If ${U}$ does not intersect ${E}$, then $\pol\notin\partial_\ss{E}$ and ${E}\setminus\{\pol\}={E}$, so ${E}\setminus\{\pol\}$ is connected.

We have
\[
{D}\;\, \subseteq \; {E}\setminus\{\pol\} \;\subseteq\;\, \rr^2\setminus\mathsf{B}_\varepsilon({A}),
\]
that is, the component ${D}$ of the subspace $\rr^2\setminus\mathsf{B}_\varepsilon({A})$ is contained in the connected subset ${E}\setminus\{\pol\}$ of this subspace.
Therefore, ${D} = {E}\setminus\{\pol\}$, and hence ${D}\subseteq{E}\subseteq{D}\cup\{\pol\}$.
Hence,
\[
\text{either } {E} = {D} \text{ or } {E} = {D} \cup \{\pol\}.
\]

Let us prove that $\partial{D}=\partial_\ss{E}$.
It suffices to show that an arbitrary point ${q}\in\ss$ either belongs or does not belong to both sets at the same time.
If ${q} = \pol$, then $q$ belongs neither to $\partial{D}$ nor to $\partial_\ss{E}$.
If ${q}\neq\pol$, then the point ${q}$ has a neighbourhood ${U}$ that does not contain $\pol$.
Since the sets ${D}$ and ${E}$ can differ only by the point $\pol$, for every neighbourhood ${V}\subseteq{U}\subseteq\ss\setminus\pol$ of the point ${q}$, we have ${V}\cap{D}={V}\cap{E}$ and ${V}\cap(\rr^2\setminus{D})={V}\cap(\ss\setminus{E})$.
Therefore, the point ${q}$ belongs to $\partial{D}$ if and only if it belongs to $\partial_\ss{E}$.

On the two-dimensional sphere $\ss$, according to Theorem~4 in~\cite[Chap.\,10, \S\,61.II, p.\,512]{kuratowski}, the following statement is true: if a locally connected continuum has no cut points, then the boundary of each component of its complement is a simple closed curve.
Thus, if we show that $\mathsf{B}_\varepsilon({A})$ is a locally connected continuum without cut points, then it follows that $\partial_\ss{E}=\partial{D}$ is a simple closed curve.
It is not difficult to show that the set ${D}$, being a component of an open subset of the plane, is open.
Then, according to Lemma~\ref{lem:complement}, ${D}$ is a component in $\rr^2\setminus\partial{D}$.
Hence, if ${D}$ is bounded, then ${D}=(\partial{D})^{-}$, and if ${D}$ is unbounded, then ${D}=(\partial{D})^{+}$.
Thus, it remains to prove that $\mathsf{B}_\varepsilon({A})$ is a locally connected continuum without cut points.

Let us show that $\mathsf{B}_\varepsilon({A})$ is a locally connected, compact set.
For every nonempty, compact set ${K}\subseteq\rr^2$ of diameter smaller than $\varepsilon$, its closed $\varepsilon$-neighbourhood $\mathsf{B}_\varepsilon({K})$ is homeomorphic~\cite[Lemma\, 1, (ii)-(iii)]{brown} to a closed disk in $\rr^2$, so it is a locally connected continuum.
Since the set ${A}$ is bounded, it can be represented as the union ${A}=\bigcup_{{i}\leqslant{n}}{A}_{i}$ of a finite number of nonempty sets of diameter less than $\varepsilon$.
For all ${i}\leqslant{n}$, the closure $\bar{A}_{i}$ is a nonempty, compact set of diameter less than $\varepsilon$.
Then $\mathsf{B}_\varepsilon(\bar{A}_{i})$ is a locally connected continuum.
Hence, $\bigcup_{{i}\leqslant{n}}\mathsf{B}_\varepsilon(\bar{A}_{i})$ is a locally connected compact \cite[Chap.\,6, \S\,49.II, p.230, Theorem\,1]{kuratowski}.
Using Remark~\ref{rmk:balls}(3) and the definition of the closed $\varepsilon$-neighbourhood we have
\[
\bigcup_{{i}\leqslant{n}}\mathsf{B}_\varepsilon(\bar{A}_{i}) = \bigcup_{{i}\leqslant{n}}\mathsf{B}_\varepsilon({A}_{i}) = \mathsf{B}_\varepsilon({A}).
\]

Let us show that $\mathsf{B}_\varepsilon({A})$ is connected and has no cut points.
To do this, it suffices to show that for any point ${r}$ in $\mathsf{B}_\varepsilon({A})$, the set $\mathsf{B}_\varepsilon({A})\setminus\{{r}\}$ is pathwise connected.
In this case, the set $\mathsf{B}_\varepsilon({A})$ is also pathwise connected.
Let ${t}$ and ${s}$ be two different points in $\mathsf{B}_\varepsilon({A})\setminus\{{r}\}$; we will find a path connecting these points in $\mathsf{B}_\varepsilon({A})\setminus\{{r}\}$.

Let ${t}'$ and ${s}'$ be points in $\bar{A}$ nearest, respectively, to ${t}$ and ${s}$.
Since ${r}$ is different from ${t}$ and ${s}$, there are points ${u}$ and ${v}$ in the segments $[{t},{t}']$ and $[{s},{s}']$, respectively, such that $\rho({u},{A})<\varepsilon$ and $\rho({v},{A})<\varepsilon$.
According to Remark~\ref{rmk:chain}(2), the open $\varepsilon$-neighbourhood $\mathsf{O}_\varepsilon({A})$ is pathwise connected, so there exists a path ${f}$ connecting the points ${u}$ and ${v}$ in $\mathsf{O}_\varepsilon({A})$.
It is easy to show that then there exists a path ${f}'$ connecting points ${t}$ and ${s}$ in $\mathsf{B}_\varepsilon({A})$ such that $\ran({f}')\setminus\{{t},{s}\}$ is contained in $\mathsf{O}_\varepsilon({A})$.

If ${r}\notin\ran({f}')$, then ${f}'$ is the path connecting ${t}$ and ${s}$ in $\mathsf{B}_\varepsilon({A})\setminus\{{r}\}$, so we are done.
If ${r}\in\ran({f}')$, then, since ${r}$ is distinct from ${t}$ and ${s}$, we have ${r}\in\ran({f}')\setminus\{{t},{s}\}\subseteq\mathsf{O}_\varepsilon({A})$.
Then there exists $\delta>0$ such that
\[
\mathsf{B}_\delta({r})\subseteq\mathsf{O}_\varepsilon({A})\setminus\{{t},{s}\}.
\]
Let $\tilde{t}$ and $\tilde{s}$ be the <<first>> and <<last>> points in the compact set $\mathsf{B}_{\delta}({r})\cap\ran({f}')$ <<on the path ${f}'$ from ${t}$ to ${s}$>>.
Then if we replace the segment of path ${f}'$ between points $\tilde{t}$ and $\tilde{s}$ with one of the arcs of the circle $\mathsf{S}_{\delta}({r})$ connecting $\tilde{t}$ and $\tilde{s}$, then we get a new path connecting ${t}$ and ${s}$, but now in $\mathsf{B}_\varepsilon({A})\setminus\{r\}$.
\end{proof}
\section{Corollaries of the theorem}
\begin{corollary}\label{cor:pzk.in.ball}
Suppose that $\varepsilon>0$ and ${A}\subseteq\rr^2$ is a nonempty, bounded, $2\varepsilon$-chained set.
Then the $\varepsilon$-boundary of ${A}$ contains a simple closed curve ${C}$ such that
\[
\mathsf{O}_\varepsilon({A})\subseteq{C}^{-}
\quad\text{ and }\quad
\mathsf{B}_\varepsilon({A})\subseteq{C}^{-}\cup{C}\,.
\]
Moreover, if the closed $\varepsilon$-neighbourhood $\mathsf{B}_\varepsilon({A})$ is simply connected, then its boundary ${E}$ is a simple closed curve and $\,\mathsf{B}_\varepsilon({A})={E}^{-}\cup{E}$.
\end{corollary}

Note that the formula $\mathsf{O}_\varepsilon({A})\subseteq{C}^{-}$ does not turn into the equality even if the set 
$B_\varepsilon(A)$ is simply connected.
For example, in the case $A = \mathsf{S}_\varepsilon({p})$, ${p}\in\rr^2$.
\begin{proof}
The number $\varepsilon$ and the set ${A}$ satisfy the conditions of Theorem~\ref{thm}.
Let ${D}$ be an unbounded component of the subspace $\rr^2\setminus\mathsf{B}_\varepsilon({A})$.
By Theorem~\ref{thm}, ${C} \coloneq \partial{D}$ is a simple closed curve, ${C}\subseteq\mathsf{S}_\varepsilon({A})$, and ${D} = {C}^+$.

Let us show that $\mathsf{O}_\varepsilon({A})\subseteq{C}^-$ and $\mathsf{B}_\varepsilon({A})\subseteq{C}^{-}\cup{C}$.
It is true that ${C}^+ = {D} \subseteq \rr^2\setminus\mathsf{B}_\varepsilon({A})$, so $\mathsf{B}_\varepsilon({A}) \cap {C}^+ = \varnothing$, hence
\[
\mathsf{O}_\varepsilon({A})\:\subseteq\:
\mathsf{B}_\varepsilon({A})\:\subseteq\:
{C}^-\!\cup{C}\:\subseteq\: 
{C}^-\!\cup\mathsf{S}_\varepsilon({A}).
\]
The only thing left to recall is that $\mathsf{O}_\varepsilon({A})\cap\mathsf{S}_\varepsilon({A})=\varnothing$.

Let $\mathsf{B}_\varepsilon({A})$ be simply connected.
We show that ${D}=\rr^2\setminus\mathsf{B}_\varepsilon({A})$.
If not, then there exists a component ${F}$ in $\rr^2\setminus\mathsf{B}_\varepsilon({A})$ different from ${D}$; in particular, ${F}\cap{D}=\varnothing$.
Then
\[
{F} \:\subseteq\: \rr^2\setminus{D} \:=\: \rr^2\setminus{C}^+ =\: {C}^-\!\cup{C}.
\]
Hence, the component ${F}$ is bounded.
According to Theorem~\ref{thm}, the boundary $\partial{F}\subseteq\mathsf{S}_\varepsilon({A})\subseteq\mathsf{B}_\varepsilon({A})$ is a simple closed curve and $(\partial{F})^-\cap\mathsf{B}_\varepsilon({A})={F}\cap\mathsf{B}_\varepsilon({A})=\varnothing$.
Thus, $\mathsf{B}_\varepsilon({A})$ contains a simple closed curve $\partial{F}$ such that $(\partial{F})^-$ is disjoint with $\mathsf{B}_\varepsilon({A})$ --- a contradiction with simple connectedness of $\mathsf{B}_\varepsilon({A})$.

Thus, ${D}$ and $\mathsf{B}_\varepsilon({A})$ are disjoint and their union equals $\rr^2$.
Consequently, $\partial{D}=\partial(\mathsf{B}_\varepsilon({A}))$, i.e., ${C}={E}$, and hence ${E}$ is a simple closed curve.
Finally, note that:
\[
\mathsf{B}_\varepsilon({A}) \:=\: \rr^2\setminus{D} \:=\: \rr^2\setminus{C}^+=\: {C}^-\!\cup{C} \:=\: {E}^-\!\cup{E}.
\]
\end{proof}

\begin{corollary}\label{cor:boundary.sep}
Suppose that ${A}\subseteq\rr^2$ is a nonempty, bounded, $2\varepsilon$-chained set, ${B}\subseteq\rr^2$ is a nonempty, $2\big(\rho({A},{B})\,{-}\,\varepsilon\big)$-chained set, and $0<\varepsilon<\rho({A},{B})$.
Then the $\varepsilon$-boundary of ${A}$ contains a simple closed curve separating ${A}$ and ${B}$.
\end{corollary}
Note that for some ${B}$ (for example, a straight line) no $\delta$-boundary of $B$ contains a simple closed curve.

\begin{proof}
Put $\delta \coloneq \rho(A,B)-\varepsilon > 0$.
Note that, for every $p\in\rr^2$, we have
\[
\rho(p,A)+\rho(p,B)\geq\rho(A,B)=\varepsilon+\delta,
\]
therefore $\mathsf{O}_\delta({B})\subseteq\rr^2\setminus\mathsf{B}_\varepsilon({A})$.
The set ${B}$ is $2\delta$-chained, hence, according to Remark~\ref{rmk:chain}(2), $\mathsf{O}_{\delta}({B})$ is pathwise connected.
Consider the component ${D}$ of the subspace $\rr^2\setminus\mathsf{B}_\varepsilon({A})$ that contains $\mathsf{O}_{\delta}({B})$.
The number $\varepsilon$ and the sets ${A}$ and ${D}$ satisfy the conditions of Theorem~\ref{thm}, so ${C} \coloneq \partial{D}$ is a simple closed curve, ${C}\subseteq\mathsf{S}_\varepsilon({A})$, and ${D} \in \{ {C}^+,{C}^- \}$.
By construction,
\[
{A}\subseteq\mathsf{O}_\varepsilon({A})
=
\mathsf{B}_\varepsilon({A})\setminus\mathsf{S}_\varepsilon({A})
\subseteq
\mathsf{B}_\varepsilon({A})\setminus{C}
\subseteq
(\rr^2\setminus{D})\setminus{C}
=
(\rr^2\setminus{C})\setminus{D} = ({C}^+\!\cup{C}^-)\setminus{D}.
\]
Thus, either ${A}\subseteq{C}^-$ or ${A}\subseteq{C}^+$, so ${C}$ separates ${A}$ and ${D}$.
Then ${C}$ separates ${A}$ and ${B}$ because ${B}\subseteq{D}$.
\end{proof}

\begin{corollary} \label{cor:comp.sep}
Suppose that $\varepsilon > 0$ and $M \subseteq \rr^2$ is a nonempty, bounded set.
Then each pair of different $2\varepsilon$-chained components of $M$ are separated by some simple closed curve contained in the $\varepsilon$-boundary of ${M}$.
\end{corollary}

\begin{proof}
Let $A$ and $B$ be two different $2\varepsilon$-chained components of $M$.
Note that $\rho(A,B) \geq 2\varepsilon$.
Then $2(\rho(A,B) - \varepsilon) \geq 2\varepsilon$, and therefore, by Remark~\ref{rmk:chain}(3), $B$ is $2(\rho(A,B)-\varepsilon)$-chained.
Thus, the number $\varepsilon$ and the sets $A$ and $B$ satisfy the conditions of Corollary~\ref{cor:boundary.sep}, so $\mathsf{S}_\varepsilon(A)$ contains a simple closed curve separating $A$ and $B$.
Finally note that $\mathsf{S}_\varepsilon(A) \subseteq \mathsf{S}_\varepsilon(M)$.
\end{proof}

From Corollaries~\ref{cor:pzk.in.ball} and \ref{cor:comp.sep} the following curious result follows, which complements Morton Brown's claim~\cite{brown} that the $\varepsilon$-boundary of a compact subset of the plane is contained in the union of a finite number of simple closed curves:

\begin{corollary} \label{cor:pzk.in.bound}
The $\varepsilon$-boundary of a nonempty, bounded subset of the plane contains a simple closed curve for all $\varepsilon > 0$.\hfill\qed
\end{corollary}

\begin{corollary}\label{cor:pzk.sep.rho/2}
Suppose that ${A},{B}\subseteq\rr^2$ are nonempty, $\rho({A},{B})$-chained sets, ${A}$ is bounded, and $\rho({A},{B})>0$.
Then the $\rho({A},{B})/2$-boundary of ${A}$ contains a simple closed curve ${C}$ that separates ${A}$ and ${B}$ and such that 
\[\rho({C},{A}\cup{B})=\rho({C},{A})=\rho({C},{B})=\rho({A},{B})/2.\]

In particular, among the simple closed curves separating ${A}$ and ${B}$, the curve ${C}$ is a curve maximally distant from ${A}\cup{B}$.
\end{corollary}

\begin{proof}
The sets $A$ and $B$ and the number $\varepsilon\coloneq \rho(A,B)/2$ satisfy the conditions of Corollary~\ref{cor:boundary.sep}.
Hence $\mathsf{S}_\varepsilon(A)$ contains a simple closed curve $C$ separating $A$ and $B$.
In particular, $\rho({C},{A}) = \varepsilon$.

Let us show that $\rho({C},{B}) \geq \varepsilon$.
If this is not true, then there are points ${p}\in{C}$, ${q}\in{B}$ and a number $\delta>0$ such that $\rho({p},{q}) < \varepsilon-\delta$.
Since ${p}\in{C}\subseteq\mathsf{S}_\varepsilon(A)$, there exists a point ${r}\in{A}$ such that $\rho({r},{p}) < \varepsilon+\delta/2$.
But then
\[
2\varepsilon = \rho({A},{B}) \leq \rho({r},{q}) \leq \rho({r},{p}) + \rho({p},{q}) < 2\varepsilon-\delta/2,
\]
a contradiction.

Now we show that $\rho({C},{B})\leq \varepsilon$.
Since $\rho({C},{A})>0$ and $\rho({C},{B})>0$, the simple closed curve ${C}$ separates $\bar{A}$ and $\bar{B}$.
The set $\bar{A}$ is compact, so there are points ${s}\in\bar{A}$ and ${t}\in\bar{B}$ such that $\rho({s},{t}) = \rho({A},{B})=2\varepsilon$.
Thus the segment $[{s},{t}]$ intersects both components of the complement of ${C}$, and so it intersects ${C}$ as well.
Let ${u}\in{C}\cap[{s},{t}]$.
Since ${u}\in{C}\subseteq\mathsf{S}_\varepsilon(A)$, then $\rho({u},{s}) \geq \varepsilon$.
Then $\rho({u},{t}) = \rho({s},{t}) - \rho({s},{u}) \leq \rho({s},{t}) - \varepsilon = \varepsilon$.
Therefore, $\rho({C},{B}) = \rho({C},\bar{B})\leq\rho({u},{t}) \leq \varepsilon$.

Thus, $\rho({C},{A}\cup{B})=\rho({C},{A})=\rho({C},{B})=\rho({A},{B})/2$.
It remains to show that
\[
\rho({C},{A}\cup{B}) = \max\{ \rho({C}',{A}\cup{B}) : {C}' \text{ is a simple closed curve separating } A \text{ and } B \}.
\]
Suppose that a simple closed curve ${C}'$ separates ${A}$ and ${B}$.
Let ${u}'\in{C}'\cap[{s},{t}]$.
Then
\[
\rho({u}',\{{s},{t}\}) \leq \rho({s},{t})/2 = \rho({A},{B})/2 = \rho({C},{A}\cup{B}).
\]
Thus,
\[
\rho({C},{A}\cup{B}) \geq \rho({u}',\{{s},{t}\})\geq \rho({C}',\overline{{A}\cup{B}}) =\rho({C}',{A}\cup{B}).
\]
\end{proof}

\begin{remark}\label{rmk:pzk.sep.max.dist}
Suppose that a simple closed curve ${C}$ separates nonempty sets ${A}$ and ${B}$ in $\rr^2$ and $0<\rho({C},{A}\cup{B})<\rho({A},{B})/2$.
Then there exists a simple closed curve ${C}'$ that separates ${A}$ and ${B}$ and such that $\rho({C}',{A}\cup{B})>\rho({C},{A}\cup{B})$.
\end{remark}

In light of Corollary~\ref{cor:pzk.sep.rho/2} and Remark~\ref{rmk:pzk.sep.max.dist}, it is interesting to note that there exist two sets, ${A}$ and ${B}$, which can be separated by some simple closed curve, but the supremum of the distances from such curves to ${A}\cup{B}$ is not attained.
See Subsection~\ref{sec:ex.rho.ab.chain} for an example.

\begin{question}\label{qst:max.dist}
Suppose that ${A},{B}\subseteq\rr^2$ are nonempty sets, $\rho({A},{B})=1$, and there exists a simple closed curve that separates ${A}$ and ${B}$.
Is there a simple closed curve ${C}$ that separates ${A}$ and ${B}$ and such that $\rho({C},{A}\cup{B})>0$?
\end{question}

\begin{proof}[Proof of Remark~\ref{rmk:pzk.sep.max.dist}.]
We may assume that ${A}\subseteq{C}^{-}$.
Put $\delta\coloneq\rho({C},{A}\cup{B})>0$.
We have $\mathsf{O}_\delta({A})=\bigcup\{\mathsf{O}_\delta({p}):{p}\in{A}\}\subseteq{C}^{-}$, so, for every ${q}\in\mathsf{O}_\delta({A})$, there is ${p}\in{A}$ such that ${q}\in\mathsf{O}_\delta({p})\subseteq\mathsf{O}_\delta({A})$.
It follows that the Lebesgue measure of each component of $\mathsf{O}_\delta({A})$ is greater than $\pi\delta^2$.
Consequently, the number of components of $\mathsf{O}_\delta({A})$ is finite; denote them by ${A}_1,\ldots,{A}_{n}$.
For every ${m}\leq{n}$, there is a path that connects in ${C}^{-}$ some point from ${A}_1$ and some point from ${A}_{m}$ because ${C}^{-}$ is pathwise connected.
Then we can find a compact set ${D}\subseteq{C}^{-}$ such that the set 
\[
{A}'\coloneq\mathsf{O}_\delta({A})\cup{D}\subseteq{C}^{-}
\]
is connected.

There exists an open ball ${O}\subseteq\rr^2$ that contains ${C}\cup{C}^{-}$.
Put ${E}\coloneq\rr^2\setminus{O}$.
We have ${E}\subseteq{C}^{+}$ and $\mathsf{O}_\delta({B})\subseteq{C}^{+}$.
The Lebesgue measure of each component of $\mathsf{O}_\delta({B})$ is greater than $\pi\delta^2$.
So, the number of components of $\mathsf{O}_\delta({B})$ that are disjoint from ${E}$ is finite.
Then we can find a compact set ${F}\subseteq{C}^{+}$ such that the set 
\[
{B}'\coloneq{E}\cup\mathsf{O}_\delta({B})\cup{F}\subseteq{C}^{+}
\]
is connected.
There exists $\varepsilon>0$ such that 
\[
\varepsilon<\rho({D}\cup{E}\cup{F},{C})\quad\text{and}\quad\varepsilon<(\rho({A},{B})/2)-\delta.
\]
Note that $\rho({A},{B})>2(\varepsilon+\delta)$.

We show that $\rho({A}',{B}')\geqslant\varepsilon$.
Let ${p}\in{A}'$ and ${q}\in{B}'$.
The segment $[{p},{q}]$ intersects ${C}$ because ${p}\in{C}^{-}$ and ${q}\in{C}^{+}$.
Let ${r}\in[{p},{q}]\cap{C}$.
If ${p}\in{D}$, then $\rho({p},{r})>\varepsilon$.
If ${q}\in{E}\cup{F}$, then $\rho({q},{r})>\varepsilon$.
In both cases, $\rho({p},{q})>\varepsilon$.
It remains to consider the case when ${p}\in\mathsf{O}_\delta({A})$ and ${q}\in\mathsf{O}_\delta({B})$.
There are $\dot{p}\in{A}$ and $\dot{q}\in{B}$ such that $\rho(\dot{p},{p})<\delta$ and $\rho({q},\dot{q})<\delta$.
By the triangle inequality, we have 
\[
\rho(\dot{p},{p})+\rho({p},{q})+\rho({q},\dot{q})\:\geq\:\rho(\dot{p},\dot{q})\:\geq\:\rho({A},{B})\:>\:2(\varepsilon+\delta),
\]
therefore
\[
\rho({p},{q})\:>\:2\varepsilon+2\delta-\rho(\dot{p},{p})-\rho({q},\dot{q})\:>\:2\varepsilon\:>\:\varepsilon.
\]

The sets ${A}',{B}'\subseteq\rr^2$ are nonempty and connected, ${A}'$ is bounded, and $\rho({A}',{B}')>0$.
Then, by Corollary~\ref{cor:pzk.sep.rho/2}, there is a simple closed curve ${C}'$ that separates ${A}'$ and ${B}'$ and such that 
\[
\rho({C}',{A}'\cup{B}')\ =\ \rho({A}',{B}')/2\ \geq\ \varepsilon/2.
\]
Then ${C}'$ that separates ${A}$ and ${B}$ because ${A}\subseteq{A}'$ and ${B}\subseteq{B}'$.

It remains to show that $\rho({C}',{A}\cup{B})>\rho({C},{A}\cup{B}).$ We shall prove 
\[
\rho({C}',{A}\cup{B})\ \geq\ \rho({C},{A}\cup{B}) + \varepsilon/2.
\] 
Let ${p}\in{C}'$ and ${q}\in{A}\cup{B}$.
First consider the case when ${q}\in{A}$.
We have $\rho({p},{q})\geq\delta$ --- otherwise, ${p}\in\mathsf{O}_\delta({A})$, and then 
\[
0\:=\:\rho({C}',\mathsf{O}_\delta({A}))\:\geq\:\rho({C}',{A}'\cup{B}')\:\geq\:\varepsilon/2
\]
because $\mathsf{O}_\delta({A})\subseteq{A}'$.
Then there is a point ${r}\in[{p},{q}]$ such that $\rho({r},{q})=\delta$.
We have ${r}\in\overline{\mathsf{O}_\delta({A})}$, so again, now using Remark~\ref{rmk:balls}(2), we have
\[
\rho({p},{r})\:\geq\:\rho({C}',\overline{\mathsf{O}_\delta({A})})\:=\:\rho({C}',\mathsf{O}_\delta({A}))\:\geq\:\rho({C}',{A}'\cup{B}')\:\geq\:\varepsilon/2.
\]
Therefore,
\[
\rho({p},{q})\:=\:\rho({p},{r})+\rho({r},{q})\:\geq\:\varepsilon/2+\delta\:=\:\rho({C},{A}\cup{B})+\varepsilon/2.
\]
The case when ${q}\in{B}$ is identical.
\end{proof}

\section{Necessity of conditions in the theorem and its corollaries}
\label{sec:exam}
What if $\varepsilon=0$? In this case the $\varepsilon$-boundary $\mathsf{S}_{\varepsilon}(A)$ of a set ${A}$ coincides with its closure.
Therefore, questions about simple closed curves lying in $\mathsf{S}_{0}(A)$ are far from the topic of this paper.
Also note that the example of Lakes of Wada~\cite{cz,wada} shows that even if ${A}$ and ${B}$ are disjoint, simply connected domains and $\rho({A},{B})=0$, there may not exist a simple closed curve separating them.

\bigskip
The following examples show the necessity of the assumptions of Theorem~\ref{thm} and its corollaries.
\subsection{Boundedness of the set ${A}$}
The boundedness condition in Theorem~\ref{thm} and Corollaries~\ref{cor:pzk.in.ball}, \ref{cor:boundary.sep}, \ref{cor:comp.sep}, and \ref{cor:pzk.in.bound} is essential.
Indeed, if the set ${A}$ is a straight line, then its $\varepsilon$-boundary does not contain a simple closed curve.
In Corollary~\ref{cor:pzk.sep.rho/2}, the boundedness condition is also essential, since two unbounded subsets of the plane cannot be separated by a simple closed curve.
\subsection{$2\varepsilon$-chainedness of the set ${A}$}
The $2\varepsilon$-chainedness condition cannot be weakened to the $(2\varepsilon+\delta)$-chainedness condition in neither Theorem~\ref{thm} nor in Corollaries~\ref{cor:pzk.in.ball} and \ref{cor:boundary.sep}, for no $\delta>0$.
Indeed, if the set ${A}$ consists of two points at distance $2\varepsilon$, then its $\varepsilon$-boundary does not contain a simple closed curve ${C}$ such that ${A}\subseteq{C}^{-}$ or ${A}\subseteq{C}^{+}$.
\subsection{Simply connectedness of the closed $\varepsilon$-neighbourhood of the set ${A}$}
The simply connectedness condition in Corollary~\ref{cor:pzk.in.ball} is essential.
Indeed, if the set ${A}$ is a circle of radius greater than $\varepsilon$, then its closed $\varepsilon$-neighbourhood is not simply connected.
And the boundary of that $\varepsilon$-neighbourhood is not a simple closed curve.
\subsection{$2(\rho({A},{B})-\varepsilon)$-chainedness of the set ${B}$}
The $2(\rho({A},{B})-\varepsilon)$-chainedness condition of the set ${B}$ in Corollary~\ref{cor:boundary.sep} cannot be relaxed to the $(2(\rho({A},{B})-\varepsilon)+\delta)$--chainedness for no $\delta>0$.
Indeed, consider a circle of radius $2\varepsilon$.
Choose points ${p}$ and ${q}$ on it such that $\rho({p},{q})=2\varepsilon$.
Let ${A}$ be the closed arc of the circle between ${p}$ and ${q}$ whose length is greater than half the length of the circle.
Let ${B}$ consists of two different points of the perpendicular bisector of the line segment $[p,q]$, which are at distance $\mathsf{min}\{\delta/2,\varepsilon\}$ from the segment $[p,q]$.
The set ${B}$ is $(2(\rho(A,B)-\varepsilon)+\delta)$--chained, but $\mathsf{S}_\varepsilon({A})$ does not contain a simple closed curve separating ${A}$ and ${B}$.
\subsection{$\rho({A},{B})$-chainedness of the sets ${A}$ and ${B}$} \label{sec:ex.rho.ab.chain}
The $\rho({A},{B})$-chainedness condition of ${A}$ in Corollary~\ref{cor:pzk.sep.rho/2} is essential.
Indeed, let ${B}$ be the closed longer arc of a circle with ends at points ${p}$ and ${q}$ such that $\rho({p},{q})$ equals the radius of the circle.
Let the point ${s}$ be the center of the circle and the point ${t}$ be the point symmetric to ${s}$ with respect to the segment $[p,q]$.
Put ${A}\coloneq\{{s},{t}\}$.
We have 
\[
\rho({p},{q})=\rho({s},{p})=\rho({s},{q})=\rho({p},{t})=\rho({q},{t})=\rho({A},{B}).
\] 
In particular, the set ${A}$ is not $\rho({A},{B})$-chained.

Suppose that a simple closed curve ${C}$ separates ${A}$ and ${B}$.
Then ${C}$ intersects the open segments $({s},{p})$ and $({p},{t})$.
Let ${u}\in{C}\cap ({s},{p})$ and ${v}\in{C}\cap({p},{t})$.
Then the simple closed curve ${B}\cup({p},{q})$ separates the points ${u}$ and ${v}$.
It follows that ${C}$ intersects $({p},{q})$ at least twice.
Therefore, 
\[
\rho({C},{A}\cup{B})\leq\rho({C},\{{p},{q}\})<\rho({p},{q})/2=\rho({A},{B})/2.
\]
In other words, there is no simple closed curve ${C}$ separating ${A}$ and ${B}$ such that $\rho({C},{A}\cup{B})=\rho({A},{B})/2$.

For every $\varepsilon>0$, there is a simple closed curve ${C}_\varepsilon$ separating ${A}$ and ${B}$ such that 
\[
\rho({C}_\varepsilon,{A}\cup{B})>(\rho({A},{B})/2)-\varepsilon.
\]
It follows that among the simple closed curves separating ${A}$ and ${B}$, there is no curve maximally distant from ${A}\cup{B}$.

If we swap the sets ${A}$ and ${B}$ in this example, we get an example showing that the $\rho({A},{B})$-chainedness condition of the set ${B}$ in Corollary~\ref{cor:pzk.sep.rho/2} is also essential.
\section{Similar questions in $\rr^3$} \label{sec:r3}

In 1976 S. Ferry showed~\cite{ferry} that the $\varepsilon$-boundary of a set ${A}\subseteq\rr^3$ is a 2-manifold for almost all $\varepsilon$.
He also proved that if ${A}$ is a finite polyhedron in $\rr^{n}$, then its $\varepsilon$-boundary is an $({n}\,{-}\,1)$-manifold for all sufficiently small values of $\varepsilon$.
In the same paper he constructed a set ${B}\subseteq\rr^3$ such that the $\varepsilon$-boundary of ${B}$ has components which are not 2-manifolds for uncountably many $\varepsilon$.
And he constructed a Cantor set in $\rr^4$ whose $\varepsilon$-boundary is not a 3-manifold for any $\varepsilon$ between 0 and 1.
\smallskip

A subset of the plane is homeomorphic to a circle if and only if it is a compact, connected, one-dimensional manifold.
Thus, in $\rr^3$, there are two different analogues of the concept of a simple closed curve, the ``spherical'' and the ``topological'':
\begin{itemize}
\item \textit{a set that is homeomorphic to the two-dimensional sphere} and
\item \textit{a compact, connected, two-dimensional manifold}.
\end{itemize}
The following example shows that in both cases the $\rr^3$ variants of Theorem~\ref{thm} and Corollaries~\ref{cor:pzk.in.ball}--\ref{cor:pzk.in.bound} do not hold.
\begin{example}
Consider two linked circles, each with a small open arc removed:
\[
{A} \coloneq \{ ({x},{y},0)\in\rr^3 : {x}^2+{y}^2=1\text{ and }{x} \geq -1+\delta \} \quad \text{and}
\]
\[
{B} \coloneq \{ ({x},0,{z})\in\rr^3 : ({x}-1)^2 + {z}^2 = 1\text{ and } x \leq 1-\delta \},
\]
where $0<\delta\ll 1$.
Let $\varepsilon \coloneq \rho({p},{q})/2$, where ${p}$ and ${q}$ are the ends of the closed arc ${A}$.

The $\varepsilon$-boundary of ${A}$ resembles the surface of a ``sausage'' bent so that its ends touch each other.
It can be shown that the $\varepsilon$-boundary of $A$ contains no nonempty, compact, connected, two-dimensional manifold.
In particular, it does not contain any set homeomorphic to the two-dimensional sphere.
The same is true about ${B}$.
\end{example}

If we modify the above example by taking ${B}$ to be a line such that $\rho({A},{B})/2=\varepsilon$, then we can see that neither the ``spherical'' nor the ``topological'' $\rr^3$ version of Corollary~\ref{cor:pzk.sep.rho/2} holds.
It is therefore interesting to consider a slightly weaker version of Corollary~\ref{cor:pzk.sep.rho/2}: 

\begin{corollary}\label{cor:pzk.sep.rho/2.modified}
Let ${A},{B}\subseteq\rr^2$ be disjoint nonempty continua.
Then there is a simple closed curve ${C}$ that separates ${A}$ and ${B}$ and such that $\rho({C},{A}\cup{B})=\rho({A},{B})/2$.
\end{corollary}

The following example shows that the ``spherical'' $\rr^3$ analogue of Corollary~\ref{cor:pzk.sep.rho/2.modified} also fails.
\begin{example}
Consider two linked circles
\[
{A}' \coloneq \{ ({x},{y},0)\in\rr^3 : {x}^2+{y}^2=1 \}
\quad \text{and} \quad
{B}' \coloneq \{ ({x},0,{z})\in\rr^3 : ({x}-1)^2 + {z}^2 = 1\}.
\]
The sets ${A}'$ and ${B}'$ cannot be separated by a set ${C}$ homeomorphic to the 2-dimensional sphere and such that $\rho({C},{A}'\cup{B}')=\rho({A}',{B}')/2$.
\end{example}

This example does not refute the ``topological'' $\rr^3$ analogue of Corollary~\ref{cor:pzk.sep.rho/2.modified}: the sets ${A}'$ and ${B}'$ can be separated by a compact, connected, two-dimensional manifold ${C}$ such that $\rho({C},{A}'\cup{B}')=\rho({A}',{B}')/2$.
So the question of the validity of the ``topological'' $\rr^3$ version of Corollary~\ref{cor:pzk.sep.rho/2.modified} remains open:

\begin{question} \label{qst:r3}
Let ${A},{B}\subseteq\rr^3$ be disjoint (simply connected) nonempty continua.
Is there a (compact, connected) two-dimensional manifold ${C}$ such that ${A}$ and ${B}$ lie in different components of $\rr^3\setminus{C}$ and $\rho({C},{A}\cup{B})=\rho({A},{B})/2$ ?
\end{question}
\printbibliography
\end{document}